\documentclass[letterpaper,11pt]{article}
\usepackage{amsmath}
\usepackage{amsthm}
\usepackage{amsfonts}
\usepackage{amssymb}
\usepackage{latexsym}
\usepackage{epsfig,graphicx,psfrag}

\newtheorem{theorem}{Theorem}

\newtheorem{lemma}[theorem]{Lemma}
\newtheorem{definition}[theorem]{Definition}
\newtheorem{proposition}[theorem]{Proposition}
\newtheorem{assumption}[theorem]{Assumption}

\numberwithin{equation}{section}
\newenvironment{Proof}{\noindent\bf{Proof.}\rm}{\hfill$\blacksquare$\bigskip}

\newcommand{\whp}{whp}

\newcommand{\eps}{\varepsilon}
\newcommand{\rmax}{r_{\rm max}}

\newcommand{\RandDist}{{\cal{F}}_{n,m,k}}
\newcommand{\RandDistSAT}{{\cal{U}}_{n,m,k}}
\newcommand{\PlantedDist}{{\cal{P}}_{n,m,k}}

\newcommand{\scO}{\mathcal{O}}

\renewcommand{\Pr}{\operatorname{Pr}}

\newcommand{\ignore}[1]{\relax}

\def\e{\varepsilon}

\def\l{\lambda}

\title{On the diameter of the set of satisfying assignments in random satisfiable $k$-CNF formulas}
\author{ Uriel Feige\thanks{The Weizmann Institute. E-mail: {\tt uriel.feige@weizmann.ac.il}.} \and
Abraham D. Flaxman\thanks{UW Institute for Health Metrics and Evaluation. E-mail: {\tt
        abie@u.washington.edu}.}  \and Dan Vilenchik\thanks{UC Berkeley. E-mail: {\tt
        vilenchi@eecs.berkeley.edu}.}}
\begin{document}
\maketitle

\begin{abstract}
%The more clauses a random formula has, the less likely it is to be
%satisfiable. Moreover,
It is known that random $k$-CNF formulas have a so-called
satisfiability threshold at a density (namely, clause-variable
ratio) of roughly $2^k\ln 2$: at densities slightly below this
threshold almost all $k$-CNF formulas are satisfiable whereas
slightly above this threshold almost no $k$-CNF formula is
satisfiable. In the current work we consider satisfiable random
formulas, and inspect another parameter -- the diameter of the
solution space (that is the maximal Hamming distance between a
pair of satisfying assignments).
%The
%conditioning on the formula being satisfiable might have different
%effects at different densities, and hence it is not clear (to the
%authors) whether the expected diameter is a monotone
%(nonincreasing) function of the density.
It was previously shown that for all densities up to a density
slightly below the satisfiability threshold the diameter is almost
surely at least roughly $n/2$ (and $n$ at much lower densities).
At densities very much higher than the satisfiability threshold,
the diameter is almost surely zero (a very dense satisfiable
formula is expected to have only one satisfying assignment). In
this paper we show that for all densities above a density that is
slightly above the satisfiability threshold (more precisely at
ratio $(1+ \eps)2^k \ln 2$, $\eps=\eps(k)$ tending to 0 as $k$
grows) the diameter is almost surely $O(k2^{-k}n)$. This shows
that a relatively small change in the density around the
satisfiability threshold (a multiplicative $(1 + \eps)$ factor),
makes a dramatic change in the diameter. This drop in the diameter
cannot be attributed to the fact that a larger fraction of the
formulas is not satisfiable (and hence have diameter 0), because
the non-satisfiable formulas are excluded from consideration by
our conditioning that the formula is satisfiable.

%It is not clear if the property of having a certain diameter is
%monotone with respect to the addition of clauses (the property of
%satisfiability clearly is monotonically decreasing), therefore it
%is a-priori not obvious why such a diminish in diameter should
%occur. In particular, the ``standard" methods to affirm a phase
%transition with respect to some property require monotonicity and
%the underlying probability space to be a product space. Both might
%not be the case in our setting (above the threshold we consider
%the uniform distribution over satisfiable $k$-CNF formulas, which
%is not known to be a product space).

\end{abstract}

\newpage

%***************************************************************************************************
%***************************************************************************************************
\section{Introduction}
%***************************************************************************************************
%***************************************************************************************************

The computational complexity of Boolean formula satisfiability has
been the focus of intensive research for decades.  Recently, a promising
approach to understanding the algorithmic difficulty of $k$-SAT has
emerged, in the form of rigorous analysis of the structural properties of
formulas drawn at random from certain distributions.  For example, a
natural distribution which has been studied extensively is the uniform
distribution over $k$-CNF formulas with exactly $m$ clauses over $n$ variables.
We denote this distribution by $\RandDist$.  Despite its
simple description, many fundamental properties of this model are yet to be
understood. For example, the computational complexity of deciding if a
random formula is satisfiable and of finding a satisfying assignment
are both major open problems \cite{Feige,Levin}.

The clause to variable ratio $m/n$ of a formula is referred to as
the {\em density} of the formula.  The random model $\RandDist$
exhibits a ``phase transition'' in satisfiability, where sparse
formulas are likely to be satisfiable whereas dense formulas are
unlikely to be satisfiable. Moreover, this phase transition
happens at a very short density interval. There exists a
satisfiability threshold $d_k=d_k(n)$ such that $k$-CNF formulas
with density $m/n > d_k$ are not satisfiable $\whp$\footnote{We
say a sequence of events holds \emph{with high probability}
($\whp$) to mean with probability tending to $1$ as $n$ tends to
infinity.}, while formulas with $m/n < d_k$ are satisfiable
$\whp$~\cite{Friedgut}. A first-moment-method calculation provides
an upper-bound of $d_k \leq 2^k\ln 2$, and the threshold is
conjectured to be within a constant distance of this upper-bound
(for all values of $k$).  A lower-bound of $2^k\ln 2 - \scO(k)$
was established rigorously using a weighted second-moment-method
in \cite{AchiPeres}.

For a satisfiable $k$-CNF formula $F$, let $\rmax(F)$ be the
maximal Hamming distance between a pair of satisfying assignments
of $F$. In this paper we study the behavior of $\rmax(F)$ as a
function of the density. Specifically, we will consider random
satisfiable formulas, and ask what the typical value of $\rmax$ is
likely to be at various densities. Observe that as one adds more
clauses to a formula, the set of satisfying assignments can only
decrease, and hence also $\rmax$ can only decrease. This indicates
that the typical value of $\rmax$ should decrease as the density
increases. However, when the formula becomes unsatisfiable, the
formula is discarded from consideration. Since the formulas of
lowest diameter (diameter 0) are those discarded from
consideration, and their proportion increases as the density
increases, this may conceivably lead to a situation in which as
the density increases the expected diameter increases rather than
decreases. In particular, there does not seem to be an a-priori
reason why the threshold for satisfiability should correspond to a
threshold behavior also with respect to the diameter of
satisfiable formulas.

Let us review what is known about $\rmax(F)$ at densities below
the satisfiability threshold. For $m/n \leq 2^{k-1} \ln 2$ we know
that all but $o(1)$-fraction of the formulas satisfy $\rmax(F)=n$
(this is because they are satisfied as NAE-$k$-SAT instances
\cite{AchiMoore}). The results in \cite{AchiRicciTers06} imply
that for $m/n = (1-\delta)2^k\ln 2$, $\delta \in (0,1/3)$, for all
but $o(1)$-fraction of satisfiable $k$-CNF formulas $\rmax(F)$ is
at least $(\frac{1}{2}-\frac{5\delta^{1/2}}{6}-\frac{2}{k})n$
(this is true for $k \geq k_0$, $k_0=k_0(\delta)$). This large
diameter is due to the existence of many small clusters of
satisfying assignments, which are ``spread" in the space of all
$2^n$ possible assignments. Physicists conjecture that this
picture persists up to the so-called condensation point at $2^k
\ln 2 - c_k$, for some constant $c_k$, at which point the number
of remaining clusters drops to polynomial and then maybe to
constant. True or not, this conjecture does not imply that
$\rmax(F)$ becomes small, because it can remain of value roughly
$n/2$ even when there are only two clusters. For densities much
higher than the satisfiability threshold (by a factor of roughly
$\log n$), the typical value of $\rmax(F)$ is~0, because such
formulas, if satisfiable, are likely to have only one satisfying
assignment (see for example~\cite{EBSPlanted} for the case of
3-CNF). This shows that the diameter of random satisfiable
formulas does undergo a phase transition as the density increases
(starting at $n$, and eventually reaching 1), but it is not clear
whether there is any density that serves as a threshold around
which there is a sharp drop in diameter.

\noindent In this paper we show that:

\begin{theorem}\label{thm:clusteringUniform} For all $k\geq 20$ and $m/n \geq (1+0.99^k)2^k\ln2$, all but a $o(1)$-fraction of satisfiable $k$-CNF formulas $F$  with $m$ clauses over $n$ variables satisfy
$$\rmax(F) \leq 50k2^{-k}n.$$
\end{theorem}
Our result proves that there occurs a transition from a typical
structure of satisfying assignments which are wide-spread in the
$n$-dimensional binary cube, to a structure where all satisfying
assignments are typically contained in a ball of small diameter.
The window in which this phase transition occurs is contained in
$[(1-\eps_1)2^k\ln 2, (1+\eps_2)2^k\ln 2]$, where both
$\eps_1,\eps_2$ tend to 0 as $k$ grows.

Here are a few interesting observations regarding this phase transition.
\begin{enumerate}
  \item The threshold phenomenon in $\rmax$ occurs at a window of densities that lies around $2^k \ln 2$, and whose width is a low-order term w.r.t. $2^k$. Since we are considering only satisfiable $k$-CNF formulas (below or above the threshold), there is no a-prior reason for this threshold to be found in the vicinity of the satisfiability threshold (as the latter is irrelevant for such formulas). Still, as our result shows, this is the case.
  \item Since we are looking at satisfiable formulas, this is not a product distribution. Therefore some methods for establishing threshold behaviors (such as \cite{Friedgut}) are not applicable.
  \item Consider the property of having a diameter of at least $r$. This is not necessarily a monotone property of the density (at least we are not aware of an easy proof that it is). Again, this shows that some approaches to prove the existence of such threshold (such as \cite{Friedgut}) may not be applicable.
  \item Typically $\rmax=n$ for $m/n < 2^k\ln2/2$. This is because
  at those ratios most formulas are satisfiable as NAE-$k$-SAT
  formulas \cite{AchiMoore} (in which case for every satisfying assignment in the NAE manner, also its complement at distance $n$ is satisfying). Numerical calculations using tools from
  statistical physics predict that at $2^k\ln k/k$ there is a phase transition
  from a typical structure of a big connected ball of satisfying assignments into many small
  balls of satisfying assignments (which are called clusters).
  Observe that  $2^k\ln k/k < 2^k\ln 2/2$ for all $k \geq 3$,
  therefore while there is a major change in the structure of the
  solution space, $\rmax$ is not affected.
\end{enumerate}

Let us briefly discuss what happens for $k <20$.  Our approach assumes
that $(2\cdot 0.99)^k$ is a low-order term compared with $2^k$. This
is however not true (or not relevant) when $k$ is small. Also, the
fact that we have a constant like 50 in the bound on $\rmax$ makes the
result trivial for small values of $k$. On the other hand, for fixed
$k$ (say $k=3$) one can numerically estimate the value of $\rmax$ (via the same methods used in the proof of Theorem \ref{thm:clusteringUniform}, just figuring out the exact numerics instead of a rigorous, less tight, estimation that we perform). For example, for $k=3$ the numerics
show that typically $\rmax < 0.2 n$ for density $m/n=7.625$ (which is $\sim 1.375\cdot 2^k\ln 2$ for $k=3$).

%\begin{figure}
%\psfrag{x}[cl][cl]{$x$}
%\psfrag{y}[cc][cc]{$f^\star(x)$}
%\psfrag{aaaa aaaa}[cl][cl]{$\epsilon = +k2^{-k}$}
%\psfrag{bbbb bbbb}[cl][cl]{$\epsilon = -k2^{-k}$}
%\begin{center}
%\includegraphics[width=4in]{f_star_k-3_eps-2^-k}
%\end{center}
%
%\caption{Plot of $f^\star(x)$ for $k=3$ and $\eps = \pm k2^{-k}$.}
%\label{fig:k=3}
%\end{figure}

\medskip

Questions regarding the structure of the solution space guided the
development of algorithms in similar contexts in the past (two
such examples are algorithms that were developed for 3CNF formulas
with a planted solution, and the intuition that served the
development of the Survey Propagation algorithm). In this paper we
limit our study to some structural properties of the solution
space and do not address algorithmic aspects, though hopefully our
new insights can serve the algorithmic perspective at some point
as well.

More precisely, while the algorithmic and structural understanding
of below-threshold random formulas and above-threshold (for
sufficiently large, yet constant, density) is rather thorough (a
short list for the below threshold regime could be
\cite{BorderFriezeUpfal,UnitClause,ChvatalReed,ClusteringPhysicists,AchiRicciTers06,AminAchi}
and
\cite{OnTheGreedy,EBSPlanted,flaxman,UniformSAT,RWalkGoldenRatio}
for the above threshold), there is no rigorous algorithmic result
for clause-variable ratio $c > 2^k\ln 2$ when $c$ is some constant
above the satisfiability threshold, but not ``sufficiently large".
(For the special case of $k=3$ there are some experimental results
in~\cite{RWalkGoldenRatio}.)

\subsection{Techniques}
One reasonable approach to prove Theorem \ref{thm:clusteringUniform}
is to consider the \textbf{uniform distribution} over satisfiable
$k$-CNF formulas with $m$ clauses over $n$ variables, and study
$\rmax(F)$ of a random instance in that distribution. Throughout
$\RandDistSAT$ denotes the uniform distribution. More specifically,
we consider a random formula $F$ from $\RandDistSAT$ and estimate
the expected number of pairs of satisfying assignments at distance
$xn$ from each other. A similar approach was used for example in
\cite{AchiPeres,ClusteringPhysicists,AchiRicciTers06} for random formulas
in the below-threshold regime.

The major additional challenge that we face in this present work is the fact that the uniform distribution $\RandDistSAT$ is not a product space, clause appearances are
dependent, and it is unclear how to quantify this dependence. On the other hand, in the below-threshold regime, since $\whp$ a random $k$-CNF formula is satisfiable, one can study random $k$-CNF formulas instead of satisfiable ones. This distribution, which we denoted above by $\RandDist$, is very ``close" to a product space (compare with the distribution where every clause is chosen independently at random with probability $p=m/\left(2^k\binom{n}{k}\right)$, which is already a product space).

One demonstration of this technical challenge is the difficulty of answering the following
question: given a fixed assignment $\psi$, what is the probability
that it satisfies a random $F$?  If $F$ is drawn from $\RandDist$
then the answer is simple, $\Pr[\psi \models F] = (1-2^{-k})^m$. If
$F$ is drawn from $\RandDistSAT$ then giving an explicit expression (as a function of $m,n,k$) for $\Pr[\psi \models F]$ is still an open question.

We will show that for $x \geq 50k2^{-k}$ the expected number of pairs of satisfying assignments at distance $xn$ from each other is much smaller than $1/n$. Since there are at most $n$ possible
ways to choose $x$, we can use the union bound to prove that $\whp$
$F$ has the desired properties (since $\RandDistSAT$ is the uniform
distribution, showing that the property holds $\whp$ translates
immediately to a deterministic statement about all but a vanishing
fraction of satisfiable formulas).

To derive our estimate on the expected number of pairs of satisfying assignments at distance $xn$ we first analyze a different distribution which is commonly called the
\textbf{planted distribution}, and we shall denote it by $\PlantedDist$.  To generate a formula according to $\PlantedDist$, fix an assignment uniformly at random, then includes $m$ clauses
uniformly at random out of $\left(2^k-1\right)\binom{n}{k}$ clauses
that are consistent with the ``planted'' assignment.

When working with $\PlantedDist$, the clauses are nearly independent
and calculation is much easier.  We then relate the planted model
and the uniform model to obtain the desired result. The idea of
translating bounds from the planted to the uniform model was used in
\cite{AminAchi,AchiRicciTers06,ClusteringPhysicists} for the
below-threshold regime, and also in \cite{Chen03,UniformSAT} but in
a different context.

The reader may wonder at this point what happens when $m/n <
(1+0.99^k)2^k\ln 2$? Do typically all satisfying assignments lie
in a low-diameter Hamming ball all the way down to the
satisfiability threshold (or even below it)? Numerical and
rigorous (tedious) calculations that we did, whose details we omit
here, suggest that Theorem \ref{thm:clusteringUniform} can be
extended (maybe with some changes in the upper bound on $\rmax$)
down to $m/n=2^k\ln 2 +O(k)$ (which is an $O(k)$-additive term
from the satisfiability threshold). This extension is done using
the same technique of going through the planted distribution.
However, when $m/n = 2^k\ln 2 +O(k)$ this technique breaks. In
Section~\ref{sec:MinimalAssignments} we discuss this issue and
suggest another technique that may prove useful when our first
technique fails. This discussion is part of a more general
discussion about the width of the window in which the phase
transition in the values of $\rmax$ occurs.

%***************************************************************************************************
%***************************************************************************************************
\section{Relating the uniform and the planted distributions }\label{sec:UnifVsPlant}
%***************************************************************************************************
%***************************************************************************************************

Let $u_x$ be a random variable counting the number of pairs of satisfying assignments at distance $xn$ from each other that a random
formula in $\RandDistSAT$ has.
%To prove Theorem \ref{thm:clusteringUniform} it suffices to show that the probability that $u_x>0$ for $x \geq 50k2^{-k}$ is $o(n^{-1})$ (and then derive Theorem \ref{thm:clusteringUniform} by using Markov's inequality and the union bound over all at most $n$ possible values of $x$).
%We start by establishing the following connection:
Let $T$ to denote the expected number of satisfying assignments that
a random formula in $\RandDistSAT$ has (that is $T=\sum_{x}E[u_x]$),
and $f_x$ a random variables which denotes the number of satisfying
assignments at distance $xn$ from the planted assignment, had $F$
belonged to $\PlantedDist$. The following proposition allows us to
upper bound $E[u_x]$ via the more accessible quantity $E[f_x]$.
\begin{proposition}\label{prop:ExpectPairsAtDist_xn} Let $F$ be a random formula sampled according to $\RandDistSAT$, then $$E[u_{x}] = T \cdot E[f_{x}]/2.$$
\end{proposition}
(A similar approach of relating the uniform and the planted distribution can be found in \cite{ClusteringPhysicists}, though in that case the uniform distribution was the non-conditioned one).\\\\
\begin{Proof}
For two satisfying assignments $\varphi_i,\varphi_j$ we use $\delta(\varphi_i,\varphi_j)$ to denote their Hamming distance.
Consider some ordering on the $2^n$ possible assignments, and let $A_i$ be an indicator variable which is 1 if $\varphi_i$ satisfies $F$. Using this terminology,
$$u_{x}=\frac{1}{2}\sum_{i,j:\delta(\varphi_i,\varphi_j)=xn}A_i\cdot A_j.$$
Linearity of expectation gives
$$E[u_{x}]=\frac{1}{2}\sum_{i,j:\delta(\varphi_i,\varphi_j)=xn}Pr[A_i \wedge  A_j]=\frac{1}{2}\sum_{\delta(\varphi_i,\varphi_j)=xn}Pr[A_i|A_j]Pr[A_j].$$
By symmetry, the latter equals
$$2^n\cdot \frac{Pr[A_j]}{2}\cdot\sum_{i:\delta(\varphi_i,\varphi_j)=xn}Pr[A_i|A_j].$$
It remains to estimate $Pr[A_i|A_j]$. Conditioning on the event
$A_j$ means conditioning on the fixed assignment $\varphi_j$ to be
satisfying. In turn, $\RandDistSAT$ conditioned on $\varphi_j$
being a satisfying assignment means that only clauses which are
satisfied by $\varphi_j$ can be included, and by symmetry, every
set of $t$ clauses satisfied by $\varphi_j$ has the same
probability of being included. Observe that for $t=m$ this is
exactly the definition of the planted distribution $\PlantedDist$.
Therefore $\sum_{i}Pr[A_i|A_j]=E[f_{x}]$, when summing over all
assignments $\varphi_i$ at distance $xn$ from $\varphi_j$.
Furthermore,  $T=\sum_jPr[A_j]$ (now we are summing over all $2^n$
assignments), and hence $Pr[A_j]=T/2^n$. Putting everything
together we derive
$$E[u_{x}] = T\cdot E[f_{x}]/2.$$
\end{Proof}

In \cite{ClusteringPhysicists} this sort of proposition was already
enough to estimate $E[u_x]$ since $T$ can be easily calculated when
$m/n$ is below the satisfiability threshold. However in
$\RandDistSAT$, $m/n$ above the satisfiability threshold, it is not
clear how to calculate $T$. The following lemma is
then useful (the proof can also be found in \cite{UniformSAT}, and
is given here for completeness).

\begin{lemma}\label{lem:ExpecRelation} Let
$W$ be the expected number of satisfying assignments of a random $\PlantedDist$ instance. Then always $T\leq W$.
\end{lemma}
\begin{Proof}
Let $t_i$ be the number of formulas on $n$ variables and $m$
clauses which have exactly $i$ satisfying assignments. Let $p_i$
be the probability that a formula with exactly $i$ satisfying
assignments is sampled from $\RandDistSAT$, and let $q_i$ be
defined similarly for $\PlantedDist$. Observe that due to
symmetry, sampling a formula from $\PlantedDist$ is exactly
equivalent to sampling a pair $(\varphi,F)$ uniformly at random
from all pairs such that $\varphi$ is an assignment and $F$ is a
formula satisfied by $\varphi$. Hence:
%For a satisfying assignment
%$\varphi$, let $\Delta_{n,m,\varphi}$ be the number of formulas on
%$n$ variables with $m$ clauses that are satisfied by $\varphi$.
%Observe that due to symmetry $\Delta_{n,m,\varphi}$ is the same
%for every $\varphi$ -- thus we omit the $\varphi$ subscript. In
%the above notation
$$p_i = \frac{t_i}{\sum_{j=1}^{2^n}t_j}, \qquad q_i = \frac{i \cdot t_i}{\sum_{i=1}^{2^n}i\cdot t_i}.$$
%\cdot\frac{i}{2^n}\cdot\frac{1}{\Delta_{n,m}}.$$
%Further observe that $$2^n\cdot \Delta_{n,m} =
%\sum_{i=1}^{2^n}i\cdot t_i.$$ This is because every formula with
%$j$ satisfying assignments is counted exactly $j$ times in the
%product $2^n\cdot \Delta_{n,m}$.
\noindent and
$$T = \sum_{i=1}^{2^n}i\cdot p_i=\frac{\sum_{i=1}^{2^n}i\cdot
t_i}{\sum_{i=1}^{2^n}t_i},$$
$$W = \sum_{i=1}^{2^n}i\cdot q_i %=\frac{\sum_{i=1}^{2^n}i^2\cdot t_i}{2^n\cdot \Delta_{n,m}}
=\frac{\sum_{i=1}^{2^n}i^2\cdot t_i}{\sum_{i=1}^{2^n}i\cdot
t_i}.$$ Therefore to prove $T \leq W$, it suffices to show
$$\left(\sum_{i=1}^{2^n}i\cdot
t_i\right)^2 \leq \left(\sum_{i=1}^{2^n}t_i\right)\cdot
\left(\sum_{i=1}^{2^n}i^2\cdot t_i\right).$$ This is just
Cauchy-Schwartz, $\left(\sum a_i\cdot b_i\right)^2\leq \left(\sum
a_i^2\right)\cdot \left(\sum b_i^2 \right)$, with $a_i =
\sqrt{t_i}$ and $b_i=i\cdot \sqrt{t_i}$.
\end{Proof}

%***************************************************************************************************
%***************************************************************************************************
\section{The Planted Setting}\label{sec:PlantedSetting}
%***************************************************************************************************
%***************************************************************************************************

In this section we analyze $W$ and $E[f_x]$. Recall that we use $W$ to denote the expected number of satisfying
assignments that a random formula in $\PlantedDist$ has, and $f_x$ counts the number of satisfying assignments at distance $xn$ from the planted assignment, had $F$ belonged to $\PlantedDist$.

\begin{figure}
\psfrag{x}[cl][cl]{$x$}
\psfrag{y}[cc][cc]{$f^\star(x)$}

%\psfrag{-0.02}[cr][cr]{$-0.02$}
%\psfrag{-0.04}[cr][cr]{$-0.04$}
%\psfrag{-0.06}[cr][cr]{$-0.06$}
%\psfrag{-0.08}[cr][cr]{$-0.08$}
%\psfrag{-0.1}[cr][cr]{$-0.10$}
%
%\psfrag{0.1}[cc][cc]{$0.1$}
%\psfrag{0.2}[cc][cc]{$0.2$}
%\psfrag{0.3}[cc][cc]{$0.3$}
%\psfrag{0.4}[cc][cc]{$0.4$}
%\psfrag{0.5}[cc][cc]{$0.5$}
%\psfrag{0.6}[cc][cc]{$0.6$}
%\psfrag{0.7}[cc][cc]{$0.7$}

\begin{center}
\includegraphics[width=4in]{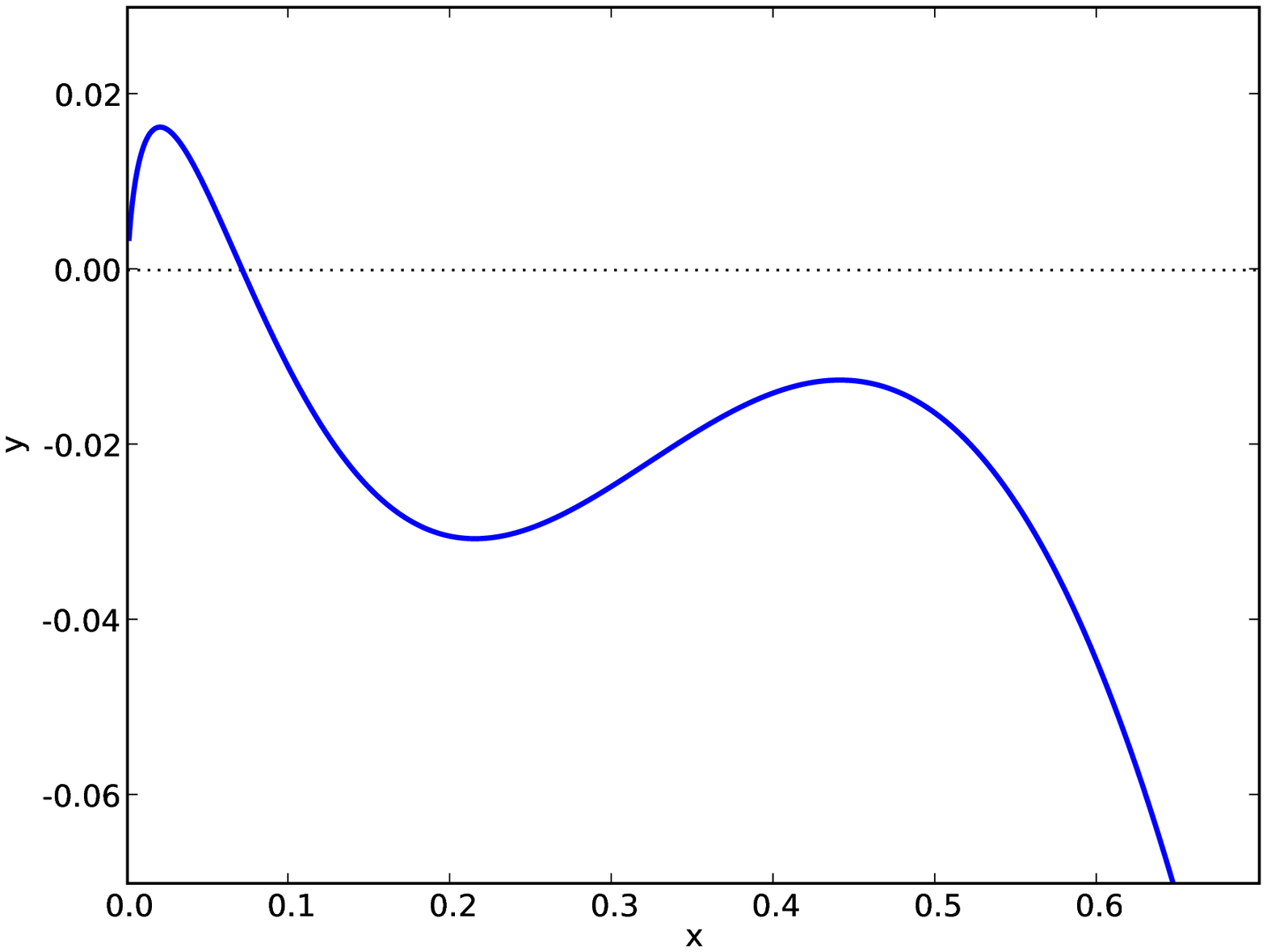}
\end{center}

\caption{Plot of $f^\star(x)$ for $k=6$ and $\eps = 2^{-k}$.}
\label{1MM-succeeds}
\end{figure}

Our analysis of $E[f_x]$ is composed of two regimes. The first is the case $x\in [0,1/k]$. In this regime we know that $E[f_x]$ changes from $\omega(1)$ to $o(1)$.  This phenomenon is depicted in Figure \ref{1MM-succeeds}. The $y$-axis in the plot is $f^\star(x)$ such that $E[f_x]=e^{f^\star(x) n}$, the $x$-axis is the Hamming distance from the planted. Therefore the transition from $E[f_x]=\omega(1)$ to $E[f_x]=o(1)$ corresponds to $f^\star(x)$ changing from positive to negative.

To translate our results to the uniform setting, it turns out that we need to have a more precise control on the rate in which $E[f_x]$ decreases once changing to $o(1)$. Therefore the analysis of that regime is more careful (Proposition \ref{prop:PlantedBreakPoint2}). Then we analyze the case $x \in [1/k,1]$. In this regime, for a suitable choice of $\eps$ (recall $m/n=(1+\eps)2^{k}\ln 2$), $E[f_x]$ is constantly $o(1)$ (in fact, exponentially small in $n$). Therefore a more crude analysis will suffice (Proposition \ref{prop:PlantedBreakPoint1}). This corresponds in Figure \ref{1MM-succeeds} to the fact that the curve is bounded away below the $x$-axis in that range.

\medskip

In this section we consider a slight modification of $\PlantedDist$. Instead of choosing $m$ clauses u.a.r., we choose $m$ clauses with repetitions. However, for $m/n=\scO(1)$, the expected number of pairs of identical clauses in $F$ (in the modified model) is $\scO(m^2/n^k)$. Thus, for $k\geq 3$ this quantity is $o(1)$. Therefore, as standard calculations show, every property that holds with probability $q$ in the modified model holds with probability $q(1+\scO(1))$ in $\PlantedDist$. Somewhat abusing notation, we will denote the modification also by $\PlantedDist$.

\medskip

\noindent Let us start with formulating $E[f_x]$ in a way which is convenient to work with.

\begin{lemma}\label{lem:FormulatingE[f_x]}
\[
E[f_x]\leq \binom{n}{xn}\cdot
\left(1 - \frac{1 - (1-x)^k}{2^k - 1}\right)^{m}.
\]
\end{lemma}
\begin{Proof}
Fix an assignment $\psi$ at distance $xn$ from the planted assignment $\varphi$. The probability that $\psi$ also satisfies $F$ can be
calculated in the following manner. Let $A$ be the set of variables on which both $\psi$ and $\varphi$ agree. $|A|=(1-x)n$.
Consider a random clause $C$ satisfied by $\varphi$; if all $k$ variables in that clause fall in $A$, then $C$ is surely satisfied by $\psi$. The probability for that is $q=\binom{(1-x)n}{k}/\binom{n}{k}$. If at least one variable falls out of $A$, which happens with probability $1-q$, then the clause is satisfied only with probability $\frac{2^k-2}{2^k-1}$. This is because there is one
way to complement the variables which is not consistent with $\psi$ but is consistent with $\varphi$. There are $\binom{n}{xn}$ ways to fix $\psi$, and therefore
$$E[f_x]=\binom{n}{xn}\left(q\cdot 1 +(1-q)\cdot\frac{2^k-2}{2^k-1}\right)^m=\binom{n}{xn}\left(\frac{2^k-2+q}{2^k-1}\right)^m
=\binom{n}{xn}\left(1-\frac{1-q}{2^k-1}\right)^m.$$
Finally, observing that $q \leq (1-x)^k$  proves the lemma.
\end{Proof}

It will be more convenient to work with the following quantity:
\begin{equation}\label{eq:DefOfFStar}
f^\star(x) \equiv \frac{\ln E[f_x]}{n}.
\end{equation}
 One can verify that
\begin{equation}\label{eq:lnOfExp}
f^\star(x)
\leq
H(x)\ln 2
+ c\ln\left(1 - \frac{1 - (1 - x)^k}{2^k - 1}\right),
\end{equation}
where $H(x)$ denotes the binary entropy measure,
\[
H(x) =
- (1-x)\log_2(1-x)
- x\log_2 x,
\]
and
$c = m/n = (1+\eps)2^{k}\ln 2$.

\medskip

To make use of Proposition \ref{prop:ExpectPairsAtDist_xn} we need to obtain tight bounds on $W$ and $E[f_x]$. In terms of $f^\star(x)$, $E[f_x]=e^{f^\star(x)n}$, therefore to prove $E[f_x]=o(1)$ it suffices to prove $f^\star(x)<0$. This is exactly what the following two propositions formally establish.

\begin{proposition}
\label{prop:PlantedBreakPoint1}
For any $k \geq 20$, $\eps \geq 0.99^k$ and $x \in [1/k,1]$,
\[
f^\star(x) \leq -50k2^{-k}
\]

\end{proposition}
\begin{Proof}
Throughout, we use the following useful upper bound on $\ln (1-x)$.
\[
\ln (1-x) \leq -x.
\]
We break the interval $[1/k,1]$ into two subintervals. Let us first consider $x \in [0.3,1]$.
Always $H(x)\ln2 \leq \ln 2$, and on the other hand, using $\log (1-x) \leq -x$,
$$c\ln\left(1 - \frac{1 - (1 - x)^k}{2^k - 1}\right) \leq -\frac{(1+\eps)2^k\ln 2}{2^k - 1}(1 - (1 - x)^k).$$
Therefore it suffices to prove that $(1+\eps)(1 - (1 - x)^k) \geq 1+\left(50k2^{-k}/\ln 2\right)$ for every $x\in [0.3,1]$.
Indeed, $$(1 - (1 - x)^k) \geq (1-0.7^k), \qquad (1+\eps)\geq (1+0.99^k).$$
One can verify that for $k \geq 20$, multiplying these two quantities is always greater than $1+\left(50k2^{-k}/\ln 2\right)$.

Let us now move the the case $x \in [1/k,0.3]$.
$H(x)$ is monotonically increasing until $x=0.5$, therefore it takes its maximal value in this interval at $x=0.3$, which gives
$H(0.3) \leq 0.266$. On the other hand $(1 - (1 - x)^k)$ takes its minimal value at $1/k$. Observe that
$(1-1/k)^{k} \leq e^{-1}$, and therefore $$(1 - (1 - x)^k) \geq 1-1/e \geq 0.6 > 0.266 > H(0.3).$$
In this case we have $f^\star(x) \leq 0.266-0.6 \leq -0.3 < 50k2^{-k}$ for every $k \geq 20$.
\end{Proof}

\begin{proposition}
\label{prop:PlantedBreakPoint2}
For any $k \geq 20$, $\e \geq 0$ and $\l \in [20,
2^k/k]$, if $x = \l 2^{-k}$ then
$
f^\star(x) \leq -\l 2^{-k}.
$
\end{proposition}
\begin{Proof}
For any $x$, we have
\[
\ln (1-x) \leq -x,
\]
and, for $0 \leq x \leq 1$,
\[
1-(1-x)^k \geq kx - \frac{k^2x^2}{2}.
\]
Thus,
\begin{align*}
&
H(x)\ln 2
+ c\ln\left(1 - \frac{1 - (1 - x)^k}{2^k - 1}\right)
\\
&\qquad=
-x\ln x - (1-x)\ln (1-x)
+ (1+\e)2^k(\ln 2)\ln\left(1-\frac{1-(1-x)^k}{2^k-1}\right)
\\
&\qquad\leq
-x\ln x + x(1-x)
- (1+\e)2^k(\ln 2)\left(\frac{1-(1-x)^k}{2^k-1}\right)
\\
&\qquad\leq
-x\ln x + x
- (1+\e)(\ln 2)\left(kx-\frac{k^2x^2}{2}\right).
\end{align*}
Substituting $\l 2^{-k}$ for $x$, this upper-bound becomes
\begin{align*}
&-x\ln x + x
- (1+\e)(\ln 2)\left(kx-\frac{k^2x^2}{2}\right)
\\
&\qquad=
\l 2^{-k} \left( k(\ln 2) - \ln\l\right)
+ \l 2^{-k}
- (1+\e)(\ln 2)\left(k\l 2^{-k} - k^2\l^2 2^{-2k-1}\right)
\\
&\qquad=
-(\l\ln \l) 2^{-k} + \l 2^{-k}
- \e(\ln 2)\left(k \l 2^{-k} - k^2 \l^2 2^{-2k-1}\right)
+ (\ln 2)k^2\l^2 2^{-2k-1}
\\
&\qquad=
-\l 2^{-k}
\left(
(\ln \l) - 1
+ \e(\ln 2)\left(k-k^2\l 2^{-k-1}\right)
- (\ln 2)k^2 \l 2^{-k-1}
\right)
\\
&\qquad=
-\l 2^{-k}
\left(
(\ln \l)\left(1 - (\ln 2)k^2\frac{\l}{\ln \l}2^{-k-1}\right) -1
+ \left(\e(\ln 2)\left(k - k^2\l2^{-k-1}\right) \right)
\right).
\end{align*}
Observe that $\l \leq 2^{k}/k$ and thus,
\[
k - k^2\l 2^{-k-1} \geq 0,
\]
and since $\eps \geq 0$ it suffices to prove that

$$(\ln \l)\left(1 - (\ln 2)k^2\frac{\l}{\ln \l}2^{-k-1}\right) -1 \geq 1.$$

Since $\l \leq 2^k/k$, and $k \geq 5$, we have
\[
(\ln 2)k^2\frac{\l}{\ln \l}2^{-k-1}
\leq
(\ln 2)k^2 \frac{2^k/k}{k(\ln 2) - \ln k}2^{-k-1}
=
(\ln 2) \frac{1}{2((\ln 2) - (\ln k)/k)}
\leq 0.65,
\]
and so it suffices to verify that
$$\ln \l \geq 2/(1-0.65),$$
which is always true for $\l \in [20,2^k/k]$
(for $k \geq 20$, $2^k/k \gg 20$).
\end{Proof}

\section{Proof of Theorem \ref{thm:clusteringUniform}}
Recall Proposition \ref{prop:ExpectPairsAtDist_xn} and  Lemma \ref{lem:ExpecRelation} which establish together
$$E[u_{x}] \leq  W \cdot E[f_{x}]/2.$$
$W$ is the expected number of satisfying assignment is the planted model, $W=\sum_{x}E[f_x]$.

\medskip

The idea of the proof is to use Propositions \ref{prop:PlantedBreakPoint1} and \ref{prop:PlantedBreakPoint2} to upper bound $W$ by looking at the largest $x$ s.t. $E[f_x]$ contributes to $W$ (that is, $E[f_x]$ is not vanishing with $n$). We shall use $x_0$ to denote this number (regardless, observe that $x_0$ is an upper bound on the diameter of the cluster region in the \emph{planted} setting). Then, to beat $W$, we take $x_1 > x_0$, so that for every $x \geq x_1$, $E[f_x]\cdot W \ll 1$. Respectively, $x_1$ uppers bounds the diameter of the cluster region in the \emph{uniform} setting. It turns out that $x_1/x_0 = \scO(k)$, and since $x_0$ scales down with $2^{-k}$, this additional factor is manageable.

\medskip

Formally, propositions \ref{prop:PlantedBreakPoint1} and \ref{prop:PlantedBreakPoint2} assert that only $x  \leq 20 \cdot 2^{-k}$ may contribute to the value of $W$. Indeed, take $x_0=20\cdot2^{-k}$, then $E[f_x] = o(n^{-1})$ for every $x \geq x_0$. For $x \leq x_0$, the total number of possible assignments (which obviously bounds the expected number of satisfying assignments) at distance $xn$ from the planted is
$$\binom{n}{xn} \leq \left(\frac{en}{xn}\right)^{xn} \leq e^{(1-\ln x)xn}.$$
This quantity is maximized for $x \leq x_0$ at $x_0$, which gives $e^{(k\ln2+1)2^{-k}n}$.  Therefore, for sufficiently large $n$,
$$W \leq o(1)+\sum_{x \leq x_0} \binom{n}{xn} \leq n e^{(k\ln 2+1)20\cdot 2^{-k}n} \leq e^{40k2^{-k}n}.$$

Now take $x_1=50k2^{-k}$ (for $k\geq 20$, $50k \leq 2^k/k$, which is the maximal $\lambda$ allowed), applying Propositions \ref{prop:PlantedBreakPoint1} and \ref{prop:PlantedBreakPoint2} once more gives that for $x \geq x_1$,
$$E[f_x] \leq e^{-50k2^{-k}n}.$$
In turn, for $x\geq x_1$
$$E[u_x] \leq W \cdot E[f_{x}]/2 \leq e^{40k2^{-k}n}\cdot e^{-50k2^{-k}n} = e^{-10k2^{-k}n}.$$
Using Markov's inequality, for $x \geq x_1$,
$$Pr[u_x >0] \leq e^{-10k2^{-k}n}.$$
Applying the union bound,
$$Pr[\exists x \geq 50k2^{-k}, u_x >0] \leq n \cdot e^{-10k2^{-k}n}=o(1).$$

%***************************************************************************************************
%***************************************************************************************************
\section{Moving even closer to the threshold}\label{sec:MinimalAssignments}
%***************************************************************************************************
%***************************************************************************************************

In the previous sections we showed that when $m/n \geq (1+0.99^k)2^k\ln2$, for $k \geq 20$, $\whp$ there are no pairs of satisfying assignments at distance greater than $50k2^{-k}$ from each other (Theorem \ref{thm:clusteringUniform}). Our approach was to consider the planted distribution and estimate $E[f_x]$ -- the expected number of satisfying assignments at distance $xn$ from the planted assignment. Then we used Proposition \ref{prop:ExpectPairsAtDist_xn} to relate this quantity to $E[u_x]$ -- the expected number of \emph{pairs} of satisfying assignments at distance $xn$ from each other (in the uniform setting). The relation we established was given (in Proposition \ref{prop:ExpectPairsAtDist_xn}) by $$E[u_x] \leq W \cdot E[f_x].$$  $W$ is the expected number of satisfying assignments in $\PlantedDist$.

\medskip

\begin{figure}
%%%% Figure 2 %%%%
\psfrag{x}[cl][cl]{$x$}
\psfrag{y}[cc][cc]{$f^\star(x)$}
%\psfrag{0.02}[cr][cr]{$0.02$}
%\psfrag{-0.02}[cr][cr]{$-0.02$}
%\psfrag{-0.04}[cr][cr]{$-0.04$}
%\psfrag{-0.06}[cr][cr]{$-0.06$}
%\psfrag{-0.08}[cr][cr]{$-0.08$}
%\psfrag{0.1}[cc][cc]{$0.1$}
%\psfrag{0.2}[cc][cc]{$0.2$}
%\psfrag{0.3}[cc][cc]{$0.3$}
%\psfrag{0.4}[cc][cc]{$0.4$}
%\psfrag{0.5}[cc][cc]{$0.5$}
%\psfrag{0.6}[cc][cc]{$0.6$}
%\psfrag{0.7}[cc][cc]{$0.7$}
\begin{center}
\includegraphics[width=4in]{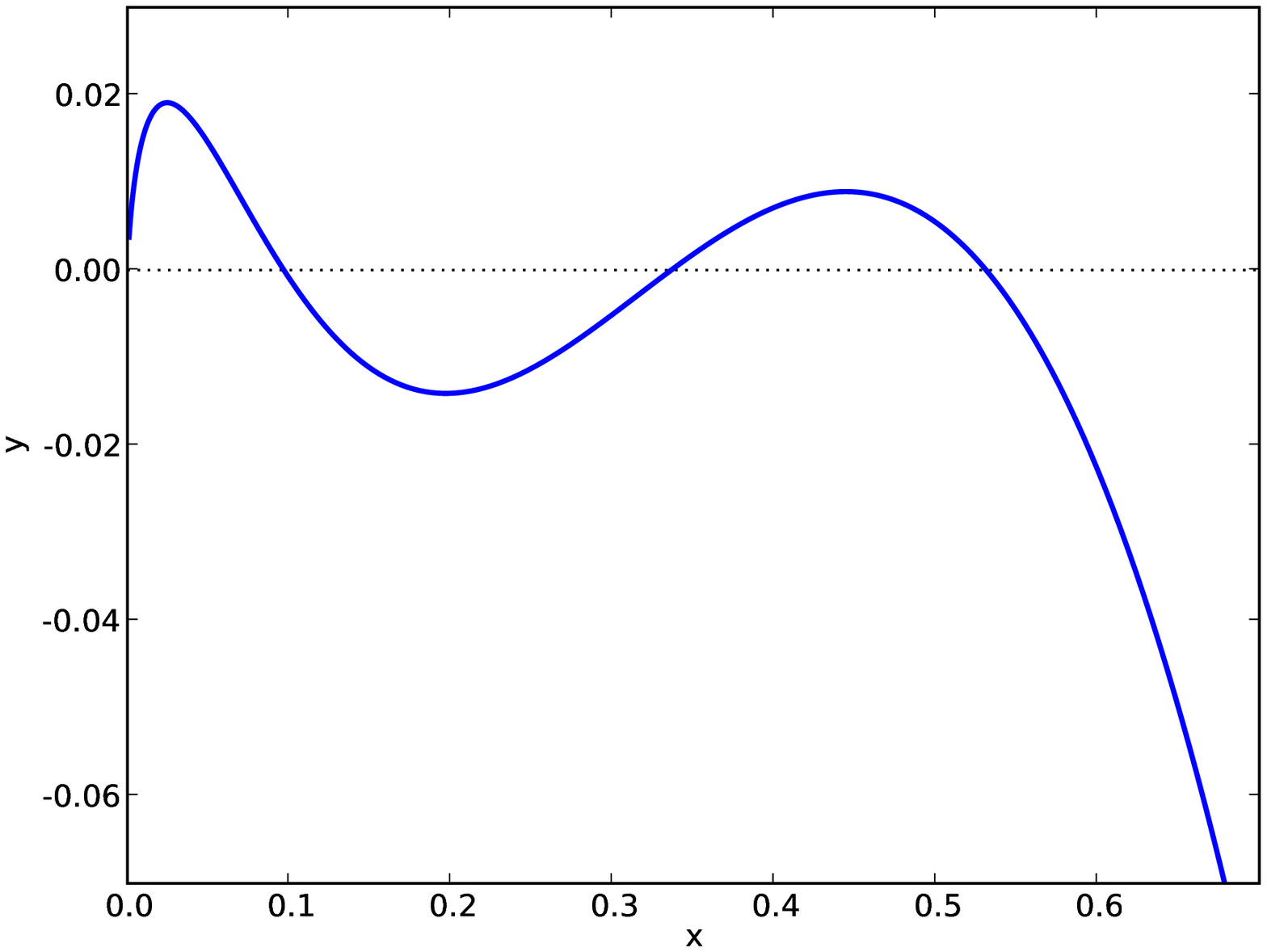}
\end{center}
\caption{Plot of $f^\star(x)$ for $k=6$ and $\eps = -2^{-k}$.}
\label{1MM-fails}

\end{figure}

Observe that $W$ is always at least 1, and therefore using this
relation to show that $E[u_x]=o(1)$ makes sense only when
$E[f_x]=o(1)$. However, using (rather tedious) calculations one
can show that when $m/n=2^k\ln 2 +O(k)$ there exists
$x\in[0.5-O(2^{-k}),0.5]$ such that $E[f_x]$ is exponentially
large in $n$ (details omitted). This phenomenon is depicted in
Figure \ref{1MM-fails} . Therefore from this density downwards our
method breaks (observe that $E[f_x]$ is monotonically decreasing
and continuous in $m/n$). This phenomenon is demonstrated in
Figures \ref{1MM-fails} vs. \ref{1MM-succeeds}.

Compare the plots in those figures. Both depict on the $y$-axis
$f^\star(x) \equiv \frac{\ln E[f_x]}{n}$, and the distance from
the planted assignment on the $x$-axis. To generate the plots we
used the estimate on $E[f_x]$ given in Lemma
\ref{lem:FormulatingE[f_x]}. Although Lemma
\ref{lem:FormulatingE[f_x]} establishes an upper bound on
$E[f_x]$, in fact for $x$ bounded away from 0 equality holds (up
to a $o(1)$ additive factor inside the parenthesis). Since
$E[f_x]$ is monotonically decreasing in $m/n$ and continuous, as
$m/n$ gets smaller, the ``hunchback" around $x=1/2$ gets closer to
the $x$-axis, and at some ratio crosses it to become positive.
This ratio occurs at $m/n =2^k\ln 2 +\scO(k)$. As $k$ grows, the
hunchback (regardless if above or below the $x$-axis) becomes
narrower, and in general is concentrated in an interval of width
$\scO(2^{-k})$ around $1/2$, with the maximum occurring at
$1/2-\scO(2^{-k})$.  We have validated these claims using a
combination of numerical and rigorous calculations (details
omitted here).

\medskip

In this section we suggest a new technique which refines the one we used. Using our refined technique we can prove for example that at some settings, even though $E[f_x]$ is exponential in $n$ (which means that our original technique fails), in fact $\whp$ $f_x=0$. Hopefully this refinement can benefit the uniform distribution as well. We do not discuss this point in the present paper.

The key to the refinement is to replace $f_x$ with another quantity which counts \emph{maximal} satisfying assignments at distance $xn$ from the planted assignment -- $f_x^{\rm max}$. This notion is similar to the notion of minimal satisfying assignments used in \cite{MinimalAssign}.

To demonstrate the power of this new technique we describe a setting where $E[f_x] \geq 1$ (which means that our original technique fails) for some $x\in[0.3,0.6]$, but $E[f_x^{\rm max}]=o(1)$ for all $x \in[0.3,0.6]$, and in that setting this will imply that $\whp$ $f_x=0$.
Formally, we prove that:

\begin{proposition}\label{prop:UntypicalBehavior} There exists a non-empty interval $(\eps_2,\eps_1)$ in which for every $\eps \in (\eps_2,\eps_1)$ and $F$ distributed according to $\PlantedDist$, $m=(1+\eps)2^k\ln 2$, there exists $x \in [0.3,0.6]$ so that $E[f_{x}]\geq 1$ while $\whp$ $f_x=0$ for every $x$ in that interval.
\end{proposition}
We choose the value $\eps_2$ carefully (we will shortly describe how), and for that $\eps_2$ we can verify numerically that
\begin{assumption} \label{assume-neg} Let $F$ be distributed according to $\PlantedDist$. If $m/n \geq (1+\eps_2)2^k\ln 2$ then $\whp$ $F$ has no satisfying assignments at distance $xn$ from the planted assignment for $0.2 \leq x \leq 0.3$ or $0.6 \leq x \leq 1.0$.
\end{assumption}
\noindent Since we are only interested in demonstrating the power of this technique, we do not care in the context of this present paper about turning it into a rigorous claim.

\medskip

Let us now formally define the notion of maximal satisfying assignments.
\begin{definition} Given a planted instance $F$ with a planted assignment $\varphi$,
we say that a satisfying assignment $\varphi'$ of $F$ is
\emph{maximal} if every assignment $\psi$ that disagrees with both $\varphi'$ and $\varphi$ on some variable $x_i$ does not satisfy $F$.
\end{definition}
In that sense $\varphi'$ is in a maximal Hamming distance from $\varphi$.
For example, if the complement of the planted also satisfies $F$, then it is maximal (in a vacant way). It is easily proven that $F$ has a satisfying assignment if and only if $F$ has a maximal satisfying assignment.

\medskip

\noindent Let $\eps_1$ be the maximal value such that for $m/n = (1+\eps_1)2^k \ln 2$ and some $x\in[0.3,0.6]$, $$E[f_{x}]\geq 1.$$ Let $\eps_2$ be the minimal value such that for $m/n =(1+\eps_2)2^k \ln 2$ and every $x\in[0.3,0.6]$ $$E[f^{\max}_{x}] \leq n^{-2}.$$

The proof of Propositions \ref{prop:PlantedBreakPoint1}  and \ref{prop:PlantedBreakPoint2} show  that $\eps_2$ always exists, and we have verified the existence of $\eps_1$ numerically. The condition $E[f^{\max}_{x}] \leq n^{-2}$ for $x \in [0.3,0.6]$ easily translates to the following claim: $\whp$ there are no maximal satisfying assignments at distance $xn$ for $x\in[0.3,0.6]$. This follows from Makrov's inequality, which gives an upper bound of $n^{-2}$ on the probability that $f^{\max}_{x}>0$ (for a fixed $x$). Now take the union bound over at most $n$ possible values of $x$.
\\\\
Before proving Proposition \ref{prop:UntypicalBehavior}, we still need to show that the interval $(\eps_2,\eps_1)$ is not empty.

\begin{proposition}\label{prop:gapInEps} $\eps_2 < \eps_1$
\end{proposition}
\begin{Proof}
Fix $x \in [0.3,0.6]$, and consider a random formula $F$ from $\PlantedDist$.
Let $M_i$ be the event that $\varphi_i$ at distance $xn$ from the planted assignment $\varphi$ is maximal, and $A_i$ the event that $\varphi_i$ satisfies $F$. Using this terminology:
\begin{align}\label{eq:ExpecMinSatAssign}
E[f^{\max}_{x}]&=\sum_{i:\delta(\varphi_i)=xn}Pr[A_i \wedge M_i]=\sum_{i}Pr[M_i | A_i]Pr[A_i]=Pr[M_i | A_i]E[f_{x}].
\end{align}
In the last step we used the fact that $Pr[M_i | A_i]$ is the same for every $\varphi_i$ by symmetry, and therefore we can pull it out in front of the summation.
It remains to estimate $Pr[M_i | A_i]$. Conditioning on the event $A_i$ in the planted model means conditioning on the fixed assignment $\varphi_i$ to be satisfying in addition to the planted assignment. In other words this means that only clauses which are satisfied by both $\varphi_i$ and $\varphi$ can be included. By symmetry, every set of $t$ clauses satisfied by both has the same probability of being included. Observe that for $t=m$ this is exactly the definition of the doubly-planted distribution (the distribution where to begin with two planted assignments are respected).

A standard approach is to consider the following variation of the doubly-planted model: pick every clause satisfied by both $\varphi_i$ and $\varphi$ w.p. $p$, where $p$ satisfies $p=m/|S|$, $S$ being the set of clauses which are satisfied by both $\varphi_i,\varphi$.
For the properties that interest us, it is straightforward to translate results between these two models.
It is also easy to see that $|S|\geq (2^k-2)\binom{n}{k}$.

Now consider a variable $s$ in $\varphi_i$ whose assignment agrees with $\varphi$, and w.l.o.g. assume it is TRUE. We call a clause $C$ \emph{$s$-qualifying} for $\varphi_i$ if it takes the form $(s \vee \ell_{y_1} \vee \ell_{y_2} \vee \ldots \vee \ell_{y_{k-1}})$,
where $\ell_{y_j}$ is a FALSE literal (over the variable $y_j$) under $\varphi_i$. If $\varphi_i$ is maximal then at least one of the  $\binom{n}{k-1}$ $s$-qualifying clauses had to be included. The probability that at least one such clause is included is at most
\begin{equation*}\label{eq:SingleVarSuppNonZero}
1-\left(1-p\right)^{\binom{n}{k-1}} \leq 1-e^{-km/(n(2^k-2))}.
\end{equation*}
Next we observe that $\varphi_i$ has at least $(1-x)n$ variables which are assigned according to $\varphi$. Also observe that the set of $s$-qualifying clauses is disjoint from
the set of $q$-qualifying clauses. Finally, for $\varphi_i$ to be maximal there must be at least one $s$-qualifying clause in $F$ for every variable $s$. The probability for that is at most
\begin{equation}\label{eq:ReductionFactor}
Pr[M_i | A_i]\leq \left(1-e^{-km/(n(2^k-2))}\right)^{(1-x)n}\leq \left(1-(1-x)e^{-km/(n(2^k-2))} \right)^n \equiv a^n,
\end{equation}
for some $a=a(k)<1$ (here we assumed that $x \in [0.3,0.6]$ and therefore $(1-x) \in [0.4,0.7]$).
Combining Equations (\ref{eq:ExpecMinSatAssign}) and (\ref{eq:ReductionFactor}) we derive

\begin{equation}\label{eq:Delta}
E[f^{\max}_{x}]\leq E[f_{x}]\cdot a^n.
\end{equation}
We claim that this implies $\eps_1 - \eps_2 \geq h$ for some $h=h(k)>0$ ($h$ actually depends on $a$, but $a$ depends only on $k$). Fix some $b=b(k)>1$ s.t. $b\cdot a < 1$ (since $a=a(k) < 1$, such $b$ exists).
Since $E[f_x]$ is continuous and decreasing in $m/n$, and by the maximality of $\eps_1$, we can find $h=h(k)>0$ s.t. $E[f_{x}]\leq b^n$ for all $x \in [0.3,0.6]$  when $m/n \leq (1+\eps_1-h)2^k\ln 2$.
On the other hand, as Equation (\ref{eq:Delta}) implies, $E[f^{\max}_{x}]\leq b^n \cdot a^n = (ab)^n \leq n^{-2}$ (for sufficiently large $n$) for all $x\in[0.3,0.6]$. By the minimality of $\eps_2$ this in particular implies that $\eps_2 \leq \eps_1-h$.
\end{Proof}\\
\begin{Proof}(Proposition \ref{prop:UntypicalBehavior})
Fix some $\eps \in (\eps_2,\eps_1)$ and consider a random formula $F$ in $\PlantedDist$ so that $m/n=(1+\eps)2^k\ln 2$. By the choice of $\eps > \eps_2$, it holds that $\whp$ $F$ has no maximal satisfying assignments at distance $xn$ from the planted assignment for $x\in [0.3,0.6]$. Assume that indeed this is the case, and also assume that Assumption \ref{assume-neg} holds.

By the choice of $\eps < \eps_1$ and the maximality of $\eps_1$, for some $x_1\in[0.3,0.6]$ indeed $E[f_{x_1}]\geq 1$. We shall now show that $f_{x}=0$ for all $x\in [0.3,0.6]$. Assume by contradiction that $f_x >0$ for some $x\in [0.3,0.6]$. Namely, there exists a satisfying assignment $\psi$ at distance $xn$ from the planted assignment, $\varphi$. Construct the assignment $\psi'$ in the following manner: while possible, flip the assignment of a variable that agrees with $\varphi$ that leaves the assignment satisfying. By construction it is clear that $\psi'$ is maximal. The crucial observation now is that at each iteration of the process we increase the distance between the current assignment and the planted by exactly one. Specifically, we start the procedure with an assignment at distance $xn$ for $x\in[0.3,0.6]$, and keep increasing the distance. If the final distance $yn$ is s.t. $y\notin [0.3,0.6]$ then at some point we've reached a satisfying assignments at distance $\geq 0.6n+1$. This
contradicts Assumption \ref{assume-neg}. Therefore we have that $\psi'$, a maximal satisfying assignment already, is at distance $yn$ for $y \in [0.3,0.6]$. This however contradicts our assumption that no maximal satisfying assignments exist at that range.
\end{Proof}

%***************************************************************************************************

\subsection*{Acknowledgements}

Large parts of this work were performed in Microsoft Research,
Redmond, Washington.

%%%%%%%%%%%%%%%%%%%%%%%%%%%%%%%%%%%%%%%%%%%%%%%%%%%%%%%%%%%%%%%%%%%%%%%%%%%%%%%%%%%%%%%%%%%%%%%%%%%%%%%%%%%%%%%%%%%%%%%%%%%%%

%%%%%%%%%%%%%%%%%%%%%%%%%%%%%%%%%%%%%%%%%%%%%%%%%%%%%%%%%%%%%%%%%%%%%%%%%%%%%%%%%%%%%%%%%%%%%%%%%%%%%%%%%%%%%%%%%%%%%%%%%%%%%


\begin{thebibliography}{10}

\bibitem{AminAchi}
D.~Achlioptas and A.~Coja-Oghlan.
\newblock Algorithmic barriers from phase transitions.
\newblock {\em preprint}.

\bibitem{AchiMoore}
D.~Achlioptas and C.~Moore.
\newblock Almost all graphs with average degree 4 are 3-colorable.
\newblock In {\em Proc. 34th ACM Symp. on Theory of Computing}, pages 199--208,
  2002.

\bibitem{AchiPeres}
D.~Achlioptas and Y.~Peres.
\newblock The threshold for random {$k$}-{SAT} is {$2^k\log 2 - O(k)$}.
\newblock {\em Journal of the AMS}, 17(4):947--973, 2004.

\bibitem{AchiRicciTers06}
D.~Achlioptas and F.~Ricci-Tersenghi.
\newblock On the solution-space geometry of random constraint satisfaction
  problems.
\newblock In {\em Proc. 38th ACM Symp. on Theory of Computing}, pages 130--139,
  2006.

\bibitem{AlonKahale97}
N.~Alon and N.~Kahale.
\newblock A spectral technique for coloring random {$3$}-colorable graphs.
\newblock {\em SIAM J. on Comput.}, 26(6):1733--1748, 1997.

\bibitem{AlonKrivSudCliqe}
N.~Alon, M.~Krivelevich, and B.~Sudakov.
\newblock Finding a large hidden clique in a random graph.
\newblock {\em Random Structures and Algorithms}, 13(3-4):457--466, 1998.

\bibitem{RWalkGoldenRatio}
E.~Amiri and Y.~Skvortsov.
\newblock Pushing random walk beyond golden ratio.
\newblock In {\em Computer Science - Theory and Applications, Second
  International Symposium on Computer Science in Russia}, pages 44--55, 2007.

\bibitem{EBSPlanted}
E.~Ben-Sasson, Y.~Bilu, and D.~Gutfreund.
\newblock Finding a randomly planted assignment in a random {$3CNF$}.
\newblock {\em manuscript}, 2002.

\bibitem{BlumSpencer}
A.~Blum and J.~Spencer.
\newblock Coloring random and semi-random {$k$}-colorable graphs.
\newblock {\em J. of Algorithms}, 19(2):204--234, 1995.

\bibitem{BorderFriezeUpfal}
A.~Z. Broder, A.~M. Frieze, and E.~Upfal.
\newblock On the satisfiability and maximum satisfiability of random 3-{CNF}
  formulas.
\newblock In {\em Proc. 4th ACM-SIAM Symp. on Discrete Algorithms}, pages
  322--330, 1993.

\bibitem{UnitClause}
M.~Chao and J.~Franco.
\newblock Probabilistic analysis of a generalization of the unit clause
  selection heuristic for the {$k$}-satisfiability problem.
\newblock {\em Information Sciences}, 51:289--314, 1990.

\bibitem{Chen03}
H.~Chen.
\newblock An algorithm for sat above the threshold.
\newblock In {\em 6th International Conference on Theory and Applications of
  Satisfiability Testing}, pages 14--24, 2003.

\bibitem{ChvatalReed}
V.~Chv{\'a}tal and B.~Reed.
\newblock Mick gets some (the odds are on his side).
\newblock In {\em Proc. 33rd IEEE Symp. on Found. of Comp. Science}, pages
  620--627, 1992.

\bibitem{UniformSAT}
A.~Coja-Oghlan, M.~Krivelevich, and D.~Vilenchik.
\newblock Why almost all satifiable $k$-{CNF} formulas are easy.
\newblock In {\em 13th conference on Analysis of Algorithms, DMTCS
  proceedings}, pages 89--102, 2007.

\bibitem{Feige}
U.~Feige.
\newblock Relations between average case complexity and approximation
  complexity.
\newblock In {\em Proc. 34th ACM Symp. on Theory of Computing}, pages 534--543,
  2002.

\bibitem{FeigeKraut}
U.~Feige and R.~Krauthgamer.
\newblock Finding and certifying a large hidden clique in a semirandom graph.
\newblock {\em Random Structures and Algorithms}, 16(2):195--208, 2000.

\bibitem{flaxman}
A.~Flaxman.
\newblock A spectral technique for random satisfiable 3{CNF} formulas.
\newblock In {\em Proc. 14th ACM-SIAM Symp. on Discrete Algorithms}, pages
  357--363, 2003.

\bibitem{Friedgut}
E.~Friedgut.
\newblock Sharp thresholds of graph properties, and the {$k$}-sat problem.
\newblock {\em J. Amer. Math. Soc.}, 12(4):1017--1054, 1999.

\bibitem{ChenFrieze}
C.~Hui and A.~M. Frieze.
\newblock Coloring bipartite hypergraphs.
\newblock In {\em Proceedings of the 5th International Conference on Integer
  Programming and Combinatorial Optimization}, pages 345--358, 1996.

\bibitem{MinimalAssign}
L.~Kirousis, E.~Kranakis, D.~Krizanc, and Y.~Stamatiou.
\newblock Approximating the unsatisfiability threshold of random formulas.
\newblock {\em Random Structures and Algorithms}, 12(3):253--269, 1998.

\bibitem{OnTheGreedy}
E.~Koutsoupias and C.~H. Papadimitriou.
\newblock On the greedy algorithm for satisfiability.
\newblock {\em Info. Process. Letters}, 43(1):53--55, 1992.

\bibitem{Levin}
L.~Levin.
\newblock Average case complete problems.
\newblock {\em SIAM J. Comput.}, 15(1):285--286, 1986.

\bibitem{ClusteringPhysicists}
M.~Mezard, T.~Mora, and R.~Zecchina.
\newblock Clustering of solutions in the random satisfiability problem.
\newblock {\em Physical Review Letters}, 94:197--205, 2005.


\end{thebibliography}
\end{document}